\documentclass[conference]{IEEEtran}

\ifCLASSINFOpdf
\else
\fi
\usepackage{times,amsmath,epsfig}
\usepackage{amssymb}
\usepackage{color,multirow}
\usepackage{booktabs}
\usepackage{breqn}
\usepackage[compress]{cite}
\usepackage{tabularx}
\usepackage{subcaption}
\usepackage{amsfonts}
\usepackage{algorithm}
\usepackage{algorithmic}
\usepackage[labelformat=simple]{subcaption}
\usepackage{setspace}
\def\hyph{--\penalty0\hskip0pt\relax}

\usepackage{tikz}
 \usepackage[usenames,dvipsnames]{pstricks}
 \usetikzlibrary{calc,trees,positioning,arrows,chains,shapes.geometric,%
    decorations.pathreplacing,decorations.pathmorphing,shapes,%
    matrix,shapes.symbols}
\usetikzlibrary{automata}
\tikzset{
>=stealth',
  punktchain/.style={
    rectangle,
    rounded corners,
    draw=black, very thick,
    text width=10em,
    minimum height=3em,
    text centered,
    on chain},
  line/.style={draw, thick, <-},
  element/.style={
    tape,
    top color=white,
    bottom color=blue!50!black!60!,
    minimum width=8em,
    draw=blue!40!black!90, very thick,
    text width=10em,
    minimum height=3.5em,
    text centered,
    on chain},
  every join/.style={->, thick,shorten >=1pt},
  decoration={brace},
  tuborg/.style={decorate},
  tubnode/.style={midway, right=2pt},
}

\definecolor {processblue}{cmyk}{0,0,0,0.17}
\definecolor{light-gray}{gray}{0.85}

\begin{document}

%
\title{On the Evaluation of Plug-in Electric Vehicle Data of a Campus Charging Network}

%
%
%

\author{%
\IEEEauthorblockN{Islam Safak Bayram \\ Qatar Environment and Energy Research Institute}
\IEEEauthorblockA{and College of Science and Engineering\\
Hamad Bin Khalifa University\\
Doha, Qatar\\
ibayram@qf.org.qa}
\and
\IEEEauthorblockN{Vahraz Zamani, Ryan Hanna, Jan Kleissl}
\IEEEauthorblockA{Department of Mechanical and Aerospace Engineering\\
University of California at San Diego, \\
La Jolla, USA\\
vazamani@ucsd.edu, rehanna@eng.ucsd.edu\\ jkleissl@eng.ucsd.edu}

}


\maketitle
\IEEEpeerreviewmaketitle

\begin{abstract}
The mass adoption of plug-in electric vehicles (PEVs) requires the deployment of public charging stations. Such facilities are expected to employ distributed generation and storage units to reduce the stress on the grid and boost sustainable transportation. While prior work has made considerable progress in deriving insights for understanding the adverse impacts of PEV chargings and how to alleviate them, a critical issue that affects the accuracy is the lack of real world PEV data. As the dynamics and pertinent design of such charging stations heavily depend on actual customer demand profile, in this paper we present and evaluate the data obtained from a $17$ node charging network equipped with Level $2$ chargers at a major North American University campus. The data is recorded for $166$ weeks starting from late $2011$. The result indicates that the majority of the customers use charging lots to extend their driving ranges. Also, the demand profile shows that there is a tremendous opportunity to employ solar generation to fuel the vehicles as there is a correlation between the peak customer demand and solar irradiation. Also, we provided a more detailed data analysis and show how to use this information in designing future sustainable charging facilities.

 \end{abstract}


\section{Introduction}
There has been a surge of interest in plug-in electric vehicles (PEVs), as electric miles are cheaper and environmentally friendly when compared to their gas-powered counterparts. This motivation is further boosted by government policies and regulations to support energy security and independence. For instance, the State of California has set a goal of putting $1.5$ million zero-emission vehicles on the road by year $2025$. Therefore, PEV sales have risen tremendously \cite{argonneEVSales}. In the last five years more than $300,000$ vehicles have been sold in the U.S. and almost $40$\% of the PEVs reside in California \cite{NSFSales,NSFSales2}. Therefore, such adoption rates have spurred the need to deploy a robust network of public charging facilities.

Typically, PEV owners use garage charging with Level $1$, $120$ AC chargers during nighttime when there is unused generation capacity. On the other hand, PEV owners living in metropolitan areas may not have access to this option. Also, the need for longer driving ranges requires the deployment of public charging facilities. During the weekdays, PEV owners can spend several hours parked at the workplace, school, or while shopping. Hence, charging service at parking lots is gaining popularity. A charging lot can be equipped with either Level $2$ single/three phase or fast DC charging technology \cite{sgc13}. On the other hand, incorporating these technologies translates to major phase, additions to the existing power grid. For instance under Level $1$ charging conditions each PEV represents a load equivalent to $50$\% of a house while under Level $2$ conditions PEV represents close to $2.5$ fold the equivalent load. 

The concern about the rapid adoption of PEVs by the utilities is that they could have disruptive impacts on the power grid. Excessive instantaneous demand can overload circuits, increase distribution transformer losses, lead to harmonics distortion, voltage deviations, and thermal loading on distribution system  \cite{eurasip,clement2010impact}. Obviously, the extent of PEV impact depends on the degree and regional density of their penetration, service requirements and the time of the day they are charged. 

For PEV adoption to occur seamlessly, the design and operation of charging facilities are of paramount importance. Corollary, there has been a growing body of literature on the design of charging facilities \cite{sgc13,power1,TSG15} and efficient load management at them \cite{lot1,lot2,jsac,TSG14}. The station designs often consider renewable generation (mostly solar power), energy storage, and efficient charge rate control. For instance, the works in \cite{jsac} and
\cite{TSG14} consider a parking lot charging facility equipped with energy storage systems to smoothen the stochastic customer demand. The studies in \cite{sgc13} and 
\cite{TSG15} provide a capacity planning framework for large-scale charging facilities. Also, utilities have started to deploy  solar-powered charging stations \cite{commercial1} and \cite{commercial2}. Moreover, the works in \cite{lot1} and \cite{lot2} propose optimization frameworks for power flow.

However, the lack of real world data limits the applicability of such studies in the US, while some work has been done in Europe \cite{intro1,intro2}. Most of the studies have to rely on simplifying assumptions, which may not be true. For instance, the common assumption on PEV service/parking time is either Gaussian or Exponential distribution, however, our results indicate that real customer behavior is far from this two assumptions. Also critical design problems such as resource provisioning and sizing of storage and PV panels heavily rely on the customer statistics, e.g., peak and spatiotemporal demand, seasonal and diurnal variations, charging duration, etc. Therefore, in this study we evaluate the data obtained from Level $2$ chargers at a major North American University campus for a period of $166$ weeks. To the best of our knowledge, this paper presents one of the most comprehensive data in this field. The main thrusts of this paper are to: (1) provide basic design and capacity planning principles for charging stations equipped with distributed generation and storage units; (2) analyze and present real-world PEV charging data on PEV charging; (3) provide guidelines for how to use the measured information to build the future sustainable charging facilities.
\begin{figure}[t]
                \centering
                \includegraphics[width=\columnwidth]{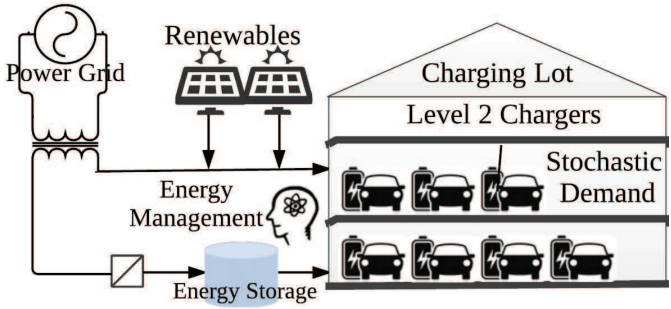}
         \caption{Illustration of PEV charging lot.}\label{model}
        \vspace{-15 pt}
\end{figure}
\section{PEV Charging Station Designs}
The design and architecture of sustainable charging facilities should carefully address the potential negative impacts of the concurrent PEV charging. The primary goal is to manage the customer demand both in time and space so that the power grid is not strained beyond its limits, and additional peaks are prevented. An overview of the systems components is presented in Fig.\ref{model}. Next, we enumerate the design parameters that enable effective demand control at PEV parking lots.

\subsection{Renewable Energy}\label{renewSection}

The recent studies conducted by Lawrence Berkeley National Laboratory and National Renewable Energy Laboratory reveal that the price of solar photovoltaic (PV) was reduced by $12$\% to $19$\% in 2013 and the trend continues to fall $3$\%-$12$\% more by the end of $2014$ in the United States \cite{vahraz1}. Moreover, Department of Energy SunShot Initiative aims to reduce PV integration cost $75$\% by 2020 \cite{vahraz2}.
As the cost of photovoltaic panels steadily drops and governments provide incentives, solar generation can be used at parking lots to fuel PEVs. The main advantage of using solar generation at Level $2$ parking lots is that customers are typically parked and stationary for 3-5 hours during their daily activities. Due to long parking durations, this method of charging has enough flexibility to exploit solar generation, e.g., the effect of intermittencies introduced by clouds, dust, etc. could be alleviated.  Also, since the electricity is generated and consumed locally, with proper control techniques it eliminates the need for transferring electricity from the far-off generator. This way the stress on the grid is also reduced. Furthermore, solar electricity is typically cheaper and cleaner than utility electricity power. For instance, according to \cite{commercial1}, the average cost per mile for solar powered vehicles is \$$0.04$, whereas the average cost for utility electricity powered and gas powered are \$$0.06$ and \$$0.20$, respectively.

\begin{figure}[t]
                \centering
                \includegraphics[width=0.7\columnwidth]{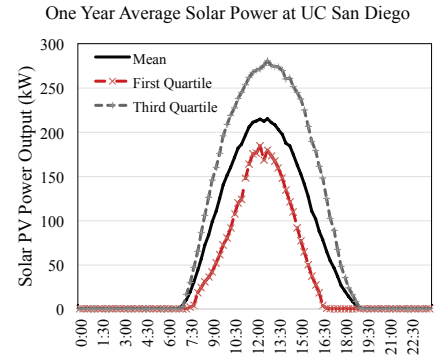}
         \caption{Solar power generation at the University campus parking deck (capacity of $300$kW).}\label{UCSDSolar}
        \vspace{-15 pt}
\end{figure}

The correlation between the seasonal customer demand and the renewable generation should be carefully addressed. For instance, one average solar PV power output of one generation facility at UC San Diego is shown in Fig. \ref{UCSDSolar}. This PV generation has $300$ kW-AC nominal capacity and is located on one of the University parking deck. As it is clear, solar generation has significant variation by presenting the first and third quartiles, i.e., $25^{th}$ percentile and $75^{th}$ percentile. This profile shows that stand-alone solar charging may not be practical. To aid the uncertainty, storage units and grid power are used.
\subsection{Energy Storage Systems (ESS)}
Energy storage units are critical components as they aid in smoothing customer demand and the integration of renewables. Moreover, due to stochasticities involved, the local generation may not perfectly align with customer demand. Storage units can act as energy buffers and decouple the time of production and demand by storing cheaper and cleaner off-peak hour electricity and delivering it during the peak load periods to handle the residual demand.

Considering their high acquisition, operation, and maintenance costs, the optimal sizing of storage units is an important design issue. Over-provisioning ESS size leads to under-utilizing costly assets and under-provisioning it taxes operation lifetime. Sizing storage units for DC fast charging stations is addressed in \cite{jsac} and \cite{TSG14}, where authors model the charging facility with a multi-dimensional Markov chain and solve the storage size problem for the peak hour. Furthermore, in \cite{lukicESS} sizing problem is solved by Monte Carlo simulations. Apparently, such studies could be improved by using data of temporal customer load profile. For instance, the sizing problem can be formulated in the following way. The station can draw a constant power to serve some base load, and the storage size could be determined to accommodate the difference between the remaining customer demand and the renewable generation. Another important aspect of storage planning is to determine the required power rating parameter, which depends on the statistics of customer charging rate, energy demand, and the length of stay. For instance, a fast-responding storage unit (e.g., supercapacitors, flywheels, etc.) may be desired to overcome the problems introduced by short-term variations in the solar generation. 

System designers should also consider the physical size of storage units as it may grow tremendously depending on the application. To illustrate this point, assume that a storage unit of size $64$kWh of Maxwell BMOD0083P048 modules (each has $26.6$Wh capacity) are chosen. Then the total volume of the ESS would be $24$m$^{3}$ for over $2400$ modules that would be required for this application. This simple example depicts how physical space can be
problematic if the station is located in regions with expensive real estate.
\subsection{Capacity Planning}
Capacity provisioning problem \hyph determining how much power to draw from the grid \hyph is an important decision problem for sustainable charging facilities. The proper resource provisioning depends on the customer demand profile and its variability. For example, if the station resides over a relatively well-confined region such as a university campus, then the variability is probably lower than a larger region such as stations located at shopping malls or interstate exchanges. Hence, resource allocation problem of the first case is relatively simpler than the second one. For instance, the work presented in \cite{TSG15} proposes a capacity planning framework for a campus size network composed of different charger technologies (Level $1$, Level $2$, DC fast etc.). They propose a capacity planning framework to guarantee a certain quality of service levels. However, since there is no measurement study, this study assumes a sinusoidal demand function for illustrative examples. Moreover, in most capacity planning studies, resource provisioning is computed respect to peak demand. On the other hand, the shape of the customer demand and the peak-to-average demand ratio can affect important design questions such as storage type, renewable portfolio, etc. Hence, it is crucial to know the daily, weekly, and the monthly demand profile. 

Another key aspect of station planning is to determine the number of physical chargers to be deployed at charging facilities. The cost per Level $2$ charger hovers between \$$600$ to \$$1000$, hence over-provisioning may be a costly option. The works presented in \cite{TSG15,TSG14,sgc13} and \cite{que1} employ queuing theory to model the charging stations. In such studies, the number of physical chargers represents ``servers'' and the operation principle investigates the case when some customers do not get the charging service immediately if all of the chargers are occupied. Hence, one of the goals is to investigate the relation between waiting time and the number of chargers. Then the design question is to find the optimal number of chargers such that customers are provided with statistical guarantees (no customer waits more than a target waiting time). Queuing-based studies heavily rely on the assumptions made on customer arrivals, departures, energy demand and the population of PEVs. In the next section, we provide detailed results of our data set.

\begin{table}[t]
\caption{Summary of PEV station design}
\small
\begin{center}
\setlength{\tabcolsep}{0.25em}
\begin{tabular}{lll}
\toprule
    Parameter & Challenges & Required Statistics \\
    \midrule
    \multicolumn{1}{c}{\multirow{2}[0]{*}{Solar }} & -Decide PV size & -Temporal customer load \\
    \multicolumn{1}{c}{} & -Forecast generation & -Charging current ranges \\ \hline
    \multicolumn{1}{c}{\multirow{2}[0]{*}{ESS}} & -Power and energy rating & -Customer demand profile \\
    \multicolumn{1}{c}{} & -Decide physical space & -PV Output \\\hline
    \multicolumn{1}{c}{\multirow{2}[0]{*}{Capacity }} & -Decide \# of Chargers & -Daily load profile\\      \multicolumn{1}{c}{} & -Capacity planning & -Arrival/Service stats.\\
   
    \bottomrule
\end{tabular}
\end{center}
\label{summary}
\vspace{-25pt}
\end{table}%

\subsection{Overview PEV Demand Control Techniques}
As the population of PEVs increase, station operators may need to apply control techniques to manage PEV demand. The customer demand can be shaped by optimally setting; (1) charging start time; (2) charging service location; and (3) charging current. An overview of techniques is presented in\cite{eurasip}. The impact and applicability of the control technique depend on customer types. For \emph{individual} customers who are interested in maximizing their benefits, incentive-based frameworks are applied to manage the spatiotemporal demand. On the other hand, companies may have \emph{PEV fleets}, who adhere to the decisions
of a \emph{central} decision maker. Such entities can up to a large extent decide charging current, service location and start time. This type of load is also considered as \emph{managable load}. For instance, the work presented in \cite{jsac} manages a hybrid population of both customer types to balance the load among neighboring stations. The work in \cite{TSG14} offers incentives to \emph{individual} customers to convince them to receive the service in a less congested station An overview is provided in Table \ref{summary}.

\section{PEV Charging Measurement}
\subsection{Charging Station Information}
The chargers are physically located in parking lots at a major North American University campus. All of the parking lots employ Level $2$, $240$ volt AC input chargers with SAE J1772 connectors. Chargers can operate up to $30$ amperes and can deliver $6.5$kW of electric power \cite{intro3}. Energy and demand metering employs ANSI 12.20 and IEC standards and the communication network uses existing cellular networks.
 
In total, there are $17$ Level $2$ chargers deployed at $9$ different parking lots. The charging nodes are located within a few miles away from each other, and exact driving distances are calculated with Google Maps and presented in Table \ref{stationInfo}. It is noteworthy that in most cases the driving distance is further than the actual physical distance. From the design point of view, this information is valuable because neighboring charging stations may be fed by the same distributed generation units (e.g., solar arrays, wind), share the same storage device, or connected to the same distribution transformer. Then the chargers of charging stations reside in small well-confined regions can be abstracted collectively as a super-station. For instance, this may apply to Stations $4-6$, as the physical distance is relatively short.

\begin{table}[t]
\caption{Driving distance (in miles) between the chargers.}
\small
\begin{center}
\setlength{\tabcolsep}{0.35em}
\begin{tabular}{cccccccccc}
\toprule
From/To& St. 1 & St. 2 & St. 3 & St. 4 & St. 5 & St. 6 & St. 7 & St. 8 & St. 9 \\
\midrule
    St. 1 & -     & 0.4   & 0.6   & 1     & 1.5   & 1     & 1.6   & 0.9   & 1.3 \\
    St. 2 & 0.4   & -     & 0.3   & 1.4   & 1.3   & 0.7   & 1.3   & 0.2   & 0.2 \\
    St. 3 & 0.6   & 0.3   & -     & 1.3   & 1.5   & 1     & 1     & 0.5   & 0.7\\
    St. 4 & 1     & 1.4   & 1.3   & -     & 0.3   & 0.3   & 2.3   & 1     & 1.4 \\
    St. 5 & 1.5   & 1.3   & 1.5   & 0.3   & -     & 0.5   & 2.5   & 1.1   & 1.6\\
    St. 6 & 1     & 0.7   & 1     & 0.3   & 0.5   & -     & 2     & 0.6   & 1\\
    St. 7 & 1.6   & 1.3   & 1     & 2.3   & 2.5   & 2     & -     & 1.4   & 1.5\\
    St. 8 & 0.9   & 0.2   & 0.5   & 1     & 1.1   & 0.6   & 1.4   & -     & 0.5\\
    St. 9 & 1.3   & 0.2   & 0.7   & 1.4   & 1.6   & 1     & 1.5   & 0.5   & - \\
    \bottomrule

\end{tabular}
\end{center}
\label{stationInfo}
\vspace{-10pt}
\end{table}%

Among these $9$ places, $2$ of them, Stations $2$ and $8$ are dedicated to University vehicles. Hence, they continue to park even after completion of the charge service. We refer to these vehicles as \emph{PEV fleet}. The remaining $7$ locations serve the \emph{public}, and we refer them as \emph{individual customers}. Detailed charging station information is presented in Table \ref{stationInfo1}.

\subsection{Data Set}
Our data set contains more than $6800$ charge events collected for $166$ weeks, starting from $14$ Nov. $2011$ to $8$ Jan. 2015.
 The collected data contains the following information \cite{etecReports}.
\begin{itemize}
\item \emph{Location name} and the \emph{address} (street name, city, state, zip code) of the charger. The chargers are at $9$ different locations.
\item \emph{Serial number} of the charger. This allows us to know the exactly which charger is in use. As a future work, this information can be used to conduct a study at a finer granularity.
\item \emph{Member guest} is the unique identification number assigned to each customer. We identify $1144$ unique visits, however, if the customer is not a member of the charger network, a new ID number will be assigned for each visit.
\item \emph{Connection and disconnection time} shows the arrival and departure of each customer. The format is in $dd/mm/yyyy$ $hh:mm:ss$.
\item \emph{Cumulative energy} shows the amount of energy (in kWh) stored at each charge service. This information is important to learn the demand statistics of vehicles and it can be used in simulation-based studies for larger scale charging networks.
\item \emph{Fee} refers the cost for the charge service in \$/kWh. The fee depends on the location and the tariffs of the utility. For instance, the rate of Level $2$ charging is $\$0.49$ per kWh in Southern California. However, the charging rate in Seattle is $10$ cents cheaper. Furthermore, there can be an additional cost in downtown areas related to parking fees.
\end{itemize}

\begin{figure*}[t]
        \centering
                \begin{subfigure}[b]{0.35\textwidth}
                \centering
                   \includegraphics[width=\columnwidth]{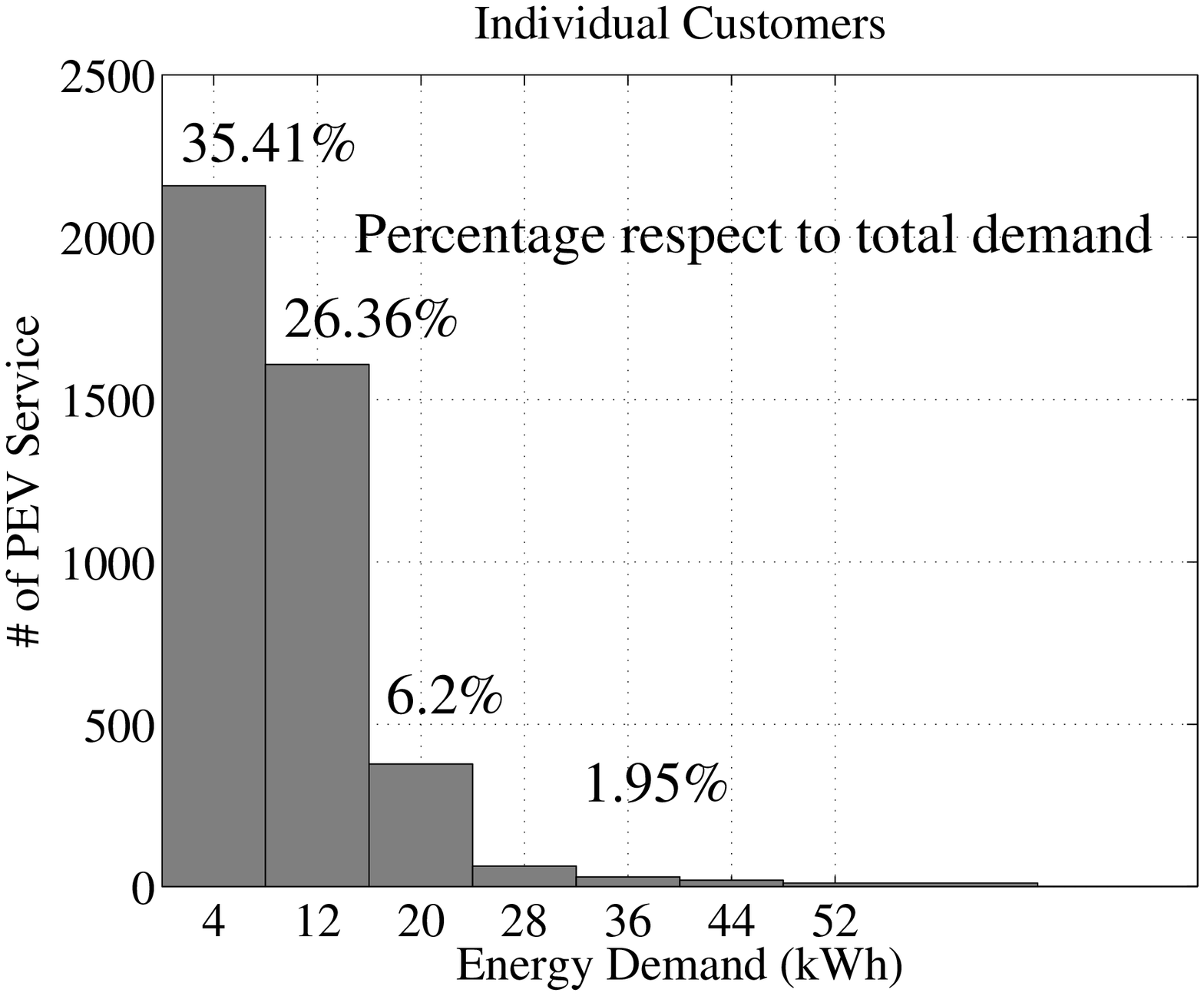}
 \caption{Energy demand of individual customers.}\label{energyIndiv}

                      \end{subfigure}
\hspace{-1.97\baselineskip}
        \begin{subfigure}[b]{0.35\textwidth}
                \centering
\includegraphics[width=\columnwidth]{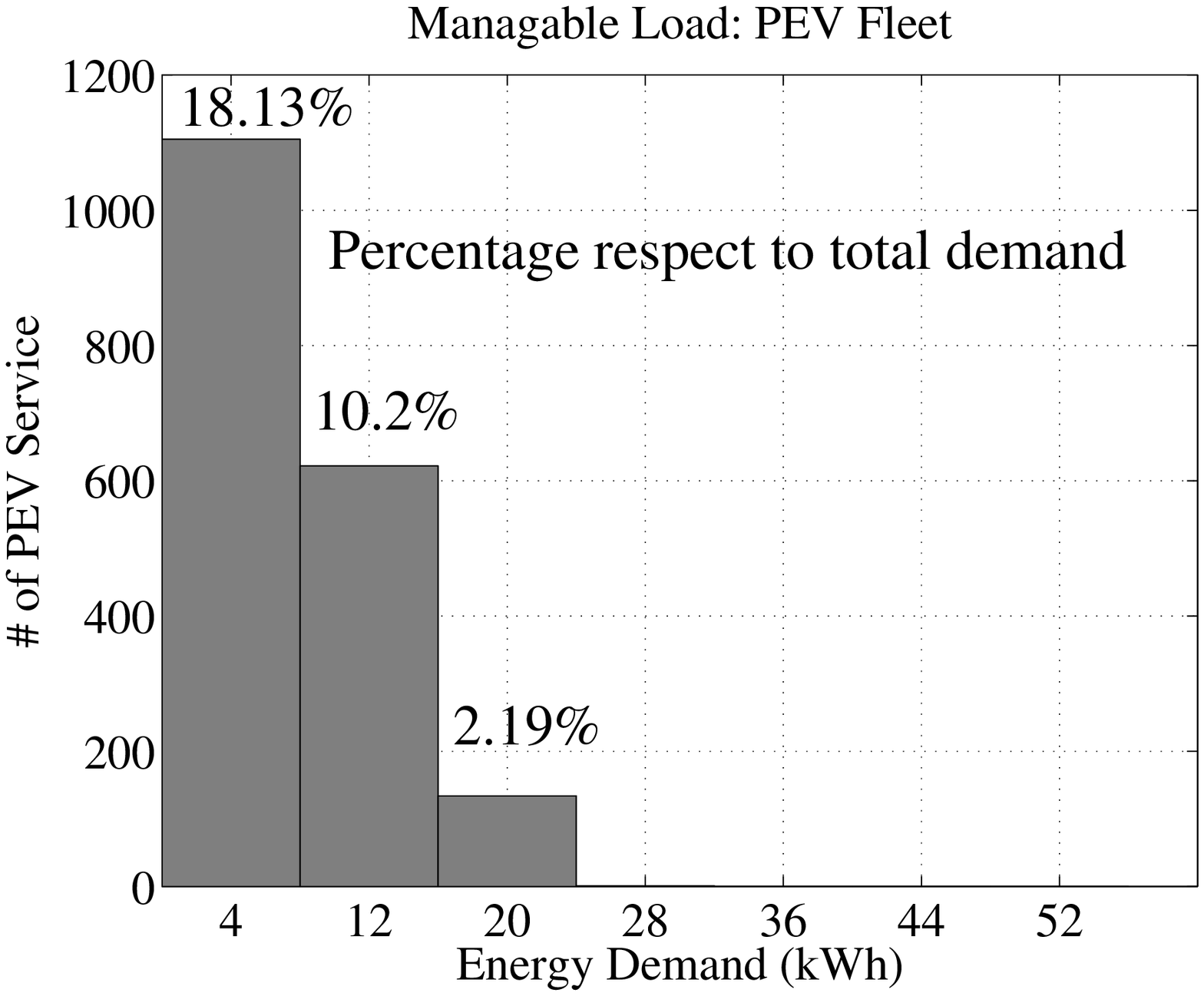}
 \caption{Energy demand of PEV fleet.}\label{energyFleet}

        \end{subfigure}%
        \hspace{-1.97\baselineskip}
        \begin{subfigure}[b]{0.35 \textwidth}
                \centering
                \includegraphics[width=\columnwidth]{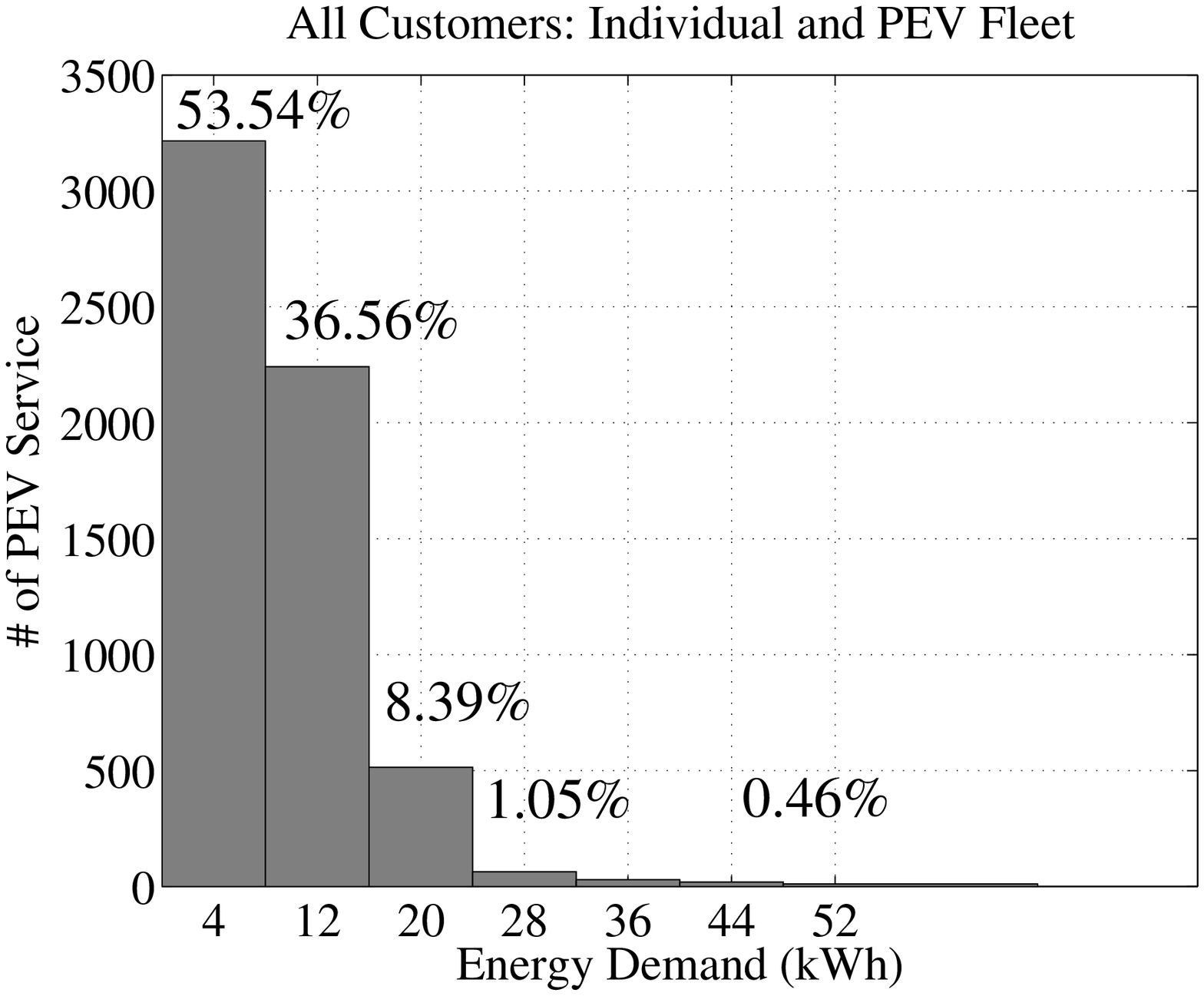}
         \caption{Aggregated energy demand.}\label{energyBoth}
        \end{subfigure}
        \caption{Energy demand distribution of single charging session. The average energy transfer is $8.53$kWh and the standard deviation is $6.49$ kWh}\label{energyDemand}
        \vspace{-6 pt}
\end{figure*}

\begin{figure*}[t]
        \centering
                \begin{subfigure}[b]{0.35\textwidth}
                \centering
                   \includegraphics[width=\columnwidth]{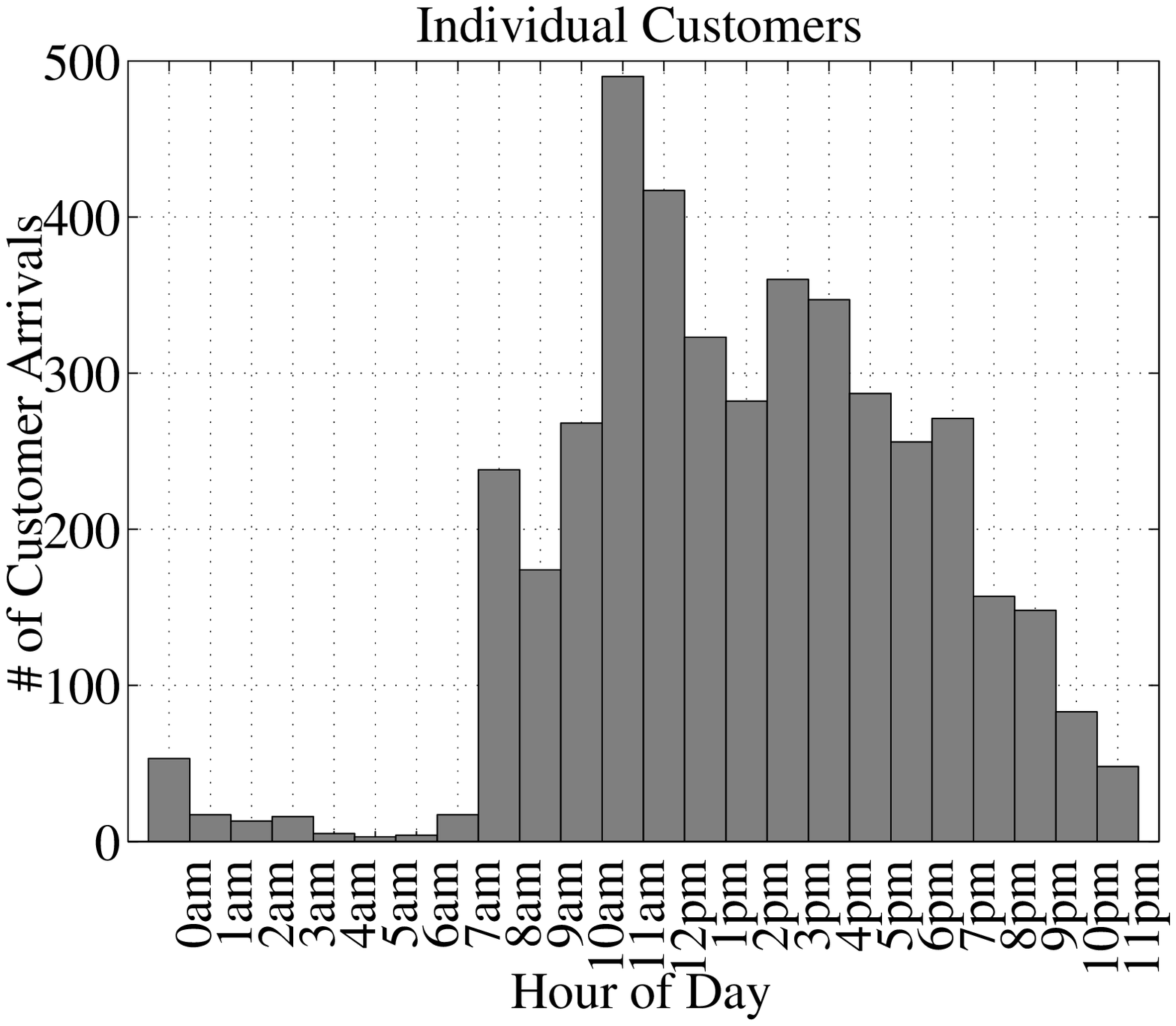}
 \caption{Hourly individual customer demand.}\label{hourIndiv}

                      \end{subfigure}
\hspace{-1.97\baselineskip}
        \begin{subfigure}[b]{0.35\textwidth}
                \centering
\includegraphics[width=\columnwidth]{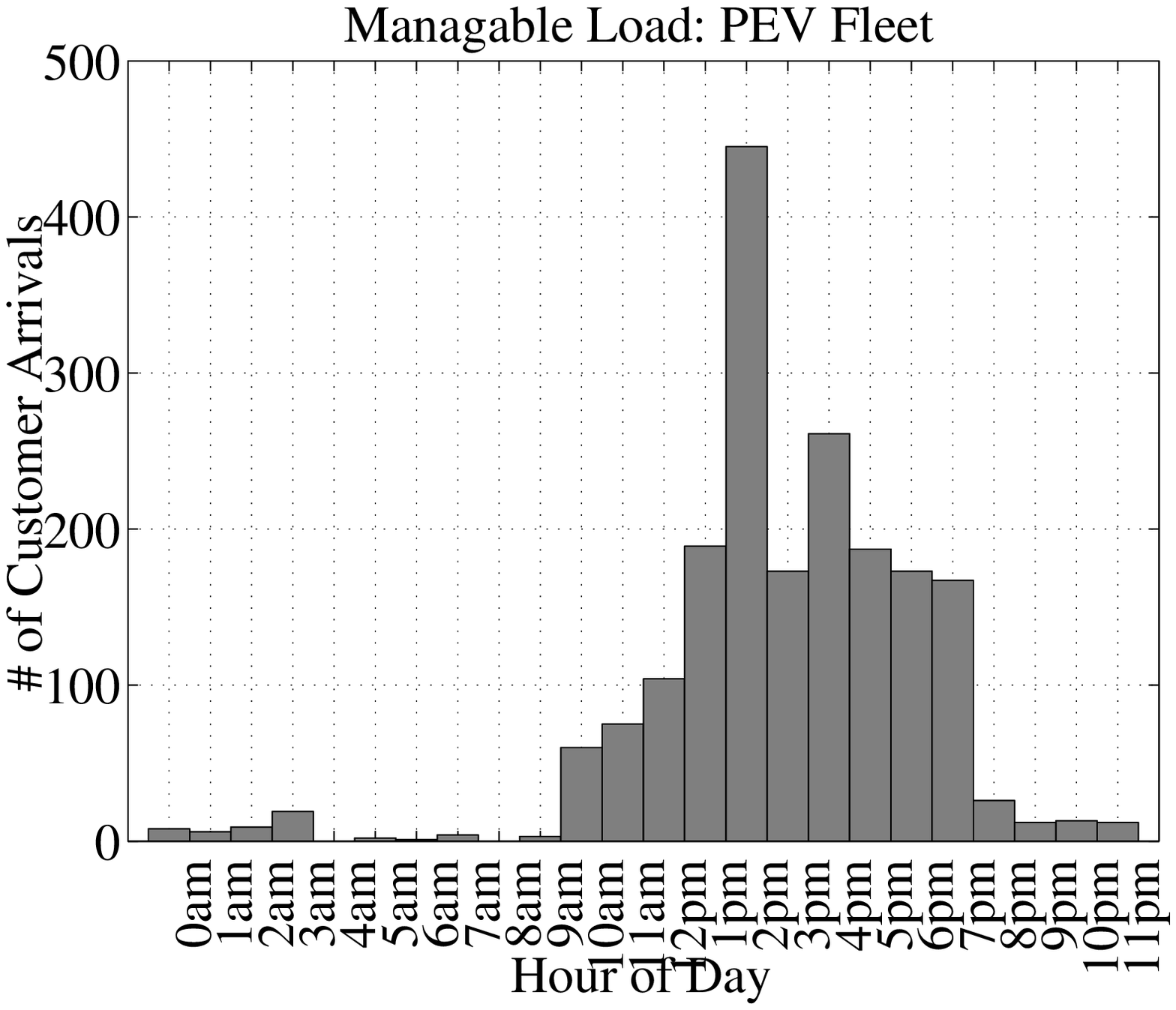}
 \caption{Hourly company PEV fleet demand.}\label{hourFleet}

        \end{subfigure}%
        \hspace{-1.97\baselineskip}
        \begin{subfigure}[b]{0.35 \textwidth}
                \centering
                \includegraphics[width=\columnwidth]{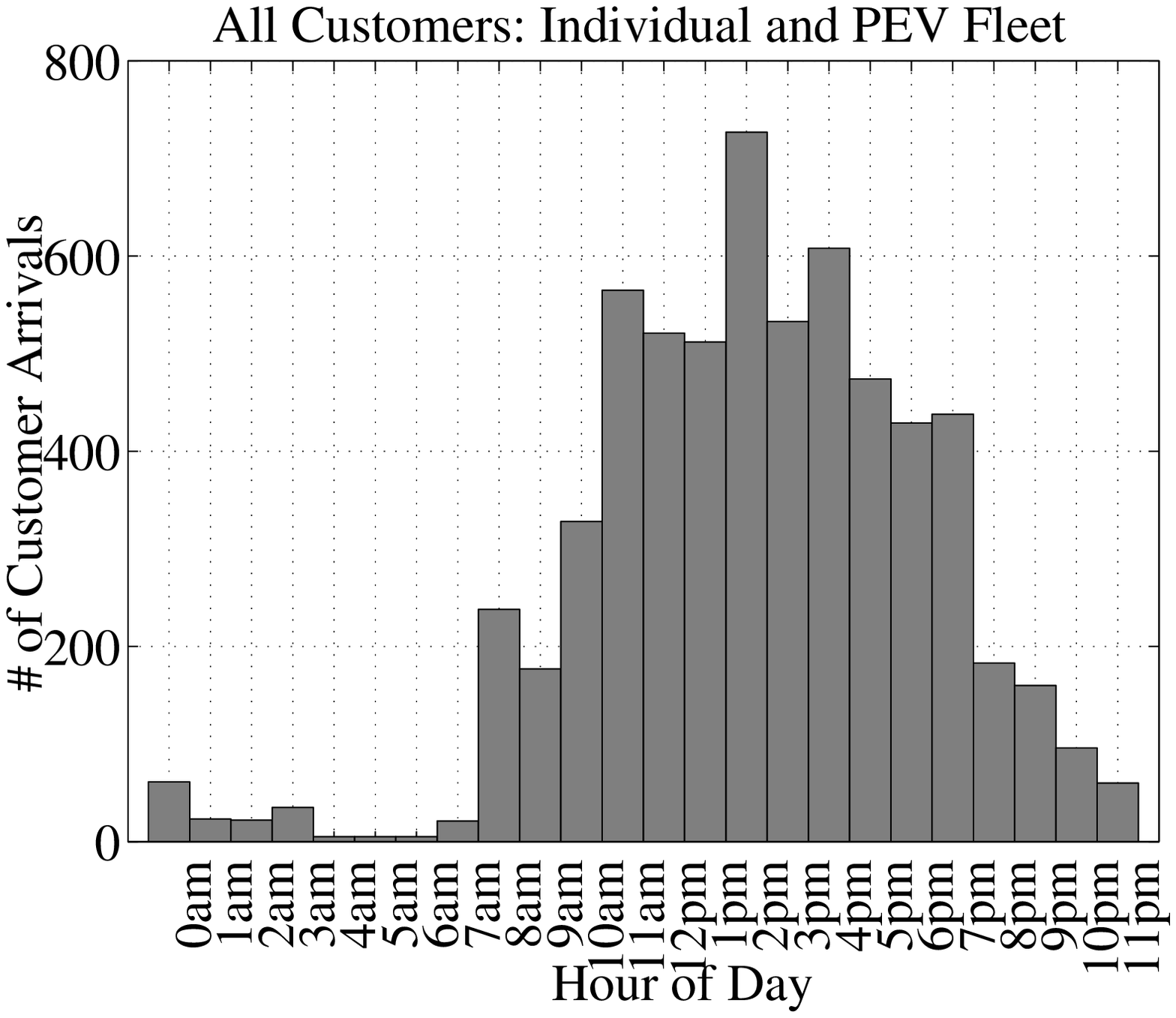}
         \caption{Hourly aggregated demand.}\label{hourBoth}
        \end{subfigure}
        \caption{PEV hourly demand}\label{hourly}
        \vspace{-9 pt}
\end{figure*}
\begin{figure*}[t]
        \centering
                \begin{subfigure}[b]{0.35\textwidth}
                \centering
                   \includegraphics[width=\columnwidth]{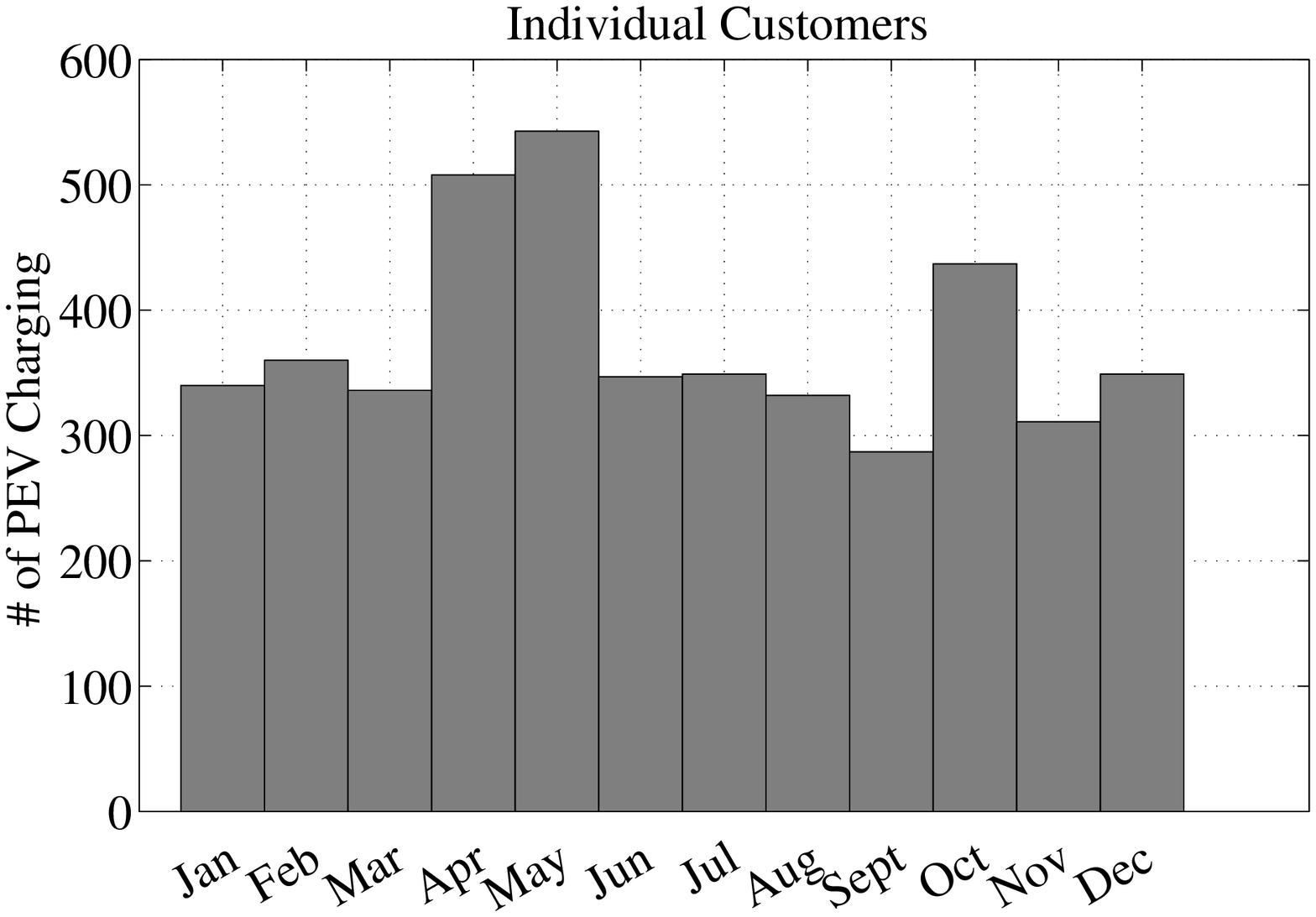}
 \caption{Monthly individual customer demand.}\label{monthIndiv}

                      \end{subfigure}
\hspace{-1.47\baselineskip}
        \begin{subfigure}[b]{0.35\textwidth}
                \centering
\includegraphics[width=\columnwidth]{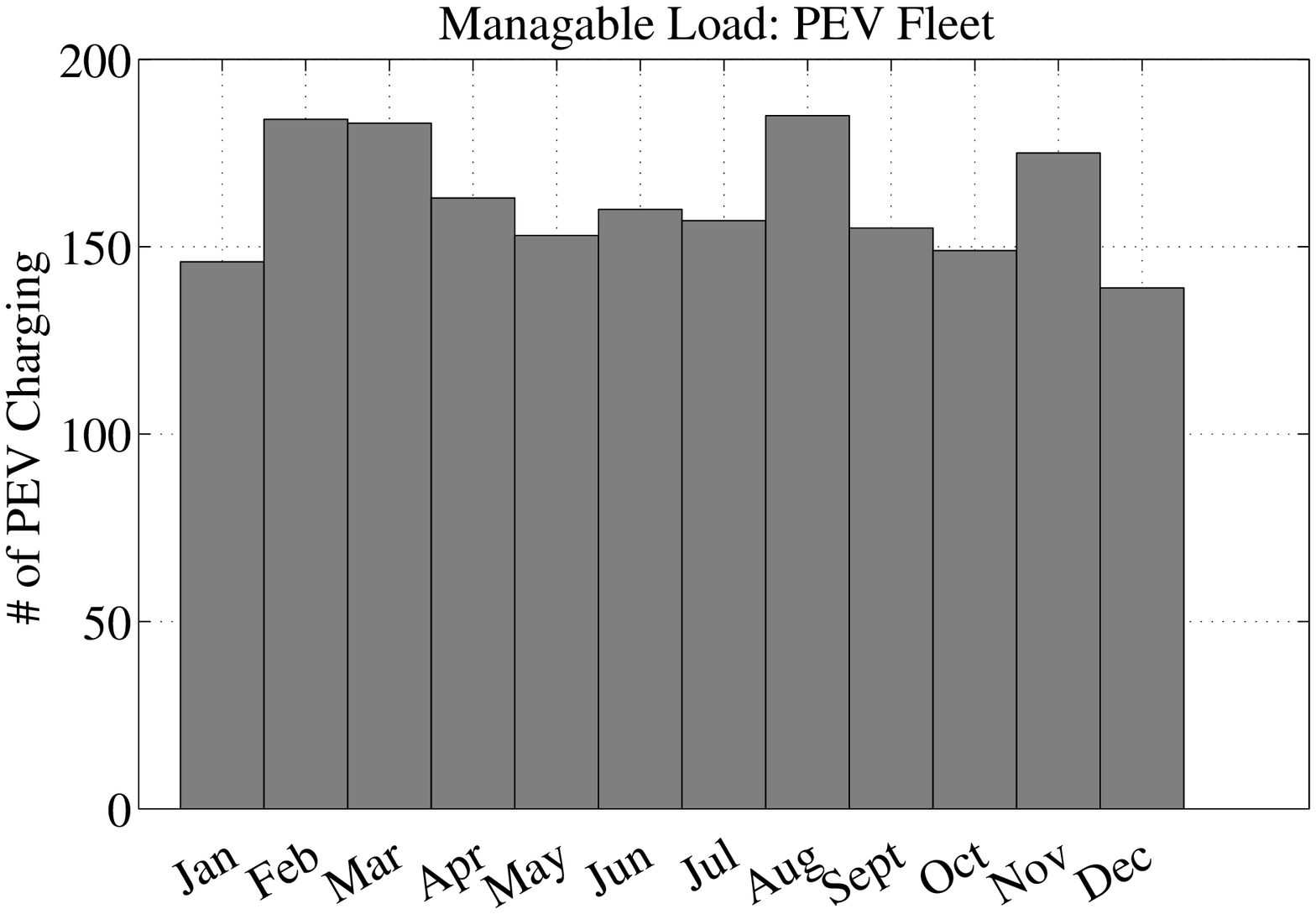}
 \caption{Monthly company PEV fleet demand.}\label{monthFleet}

        \end{subfigure}%
        \hspace{-1.47\baselineskip}
        \begin{subfigure}[b]{0.35 \textwidth}
                \centering
                \includegraphics[width=\columnwidth]{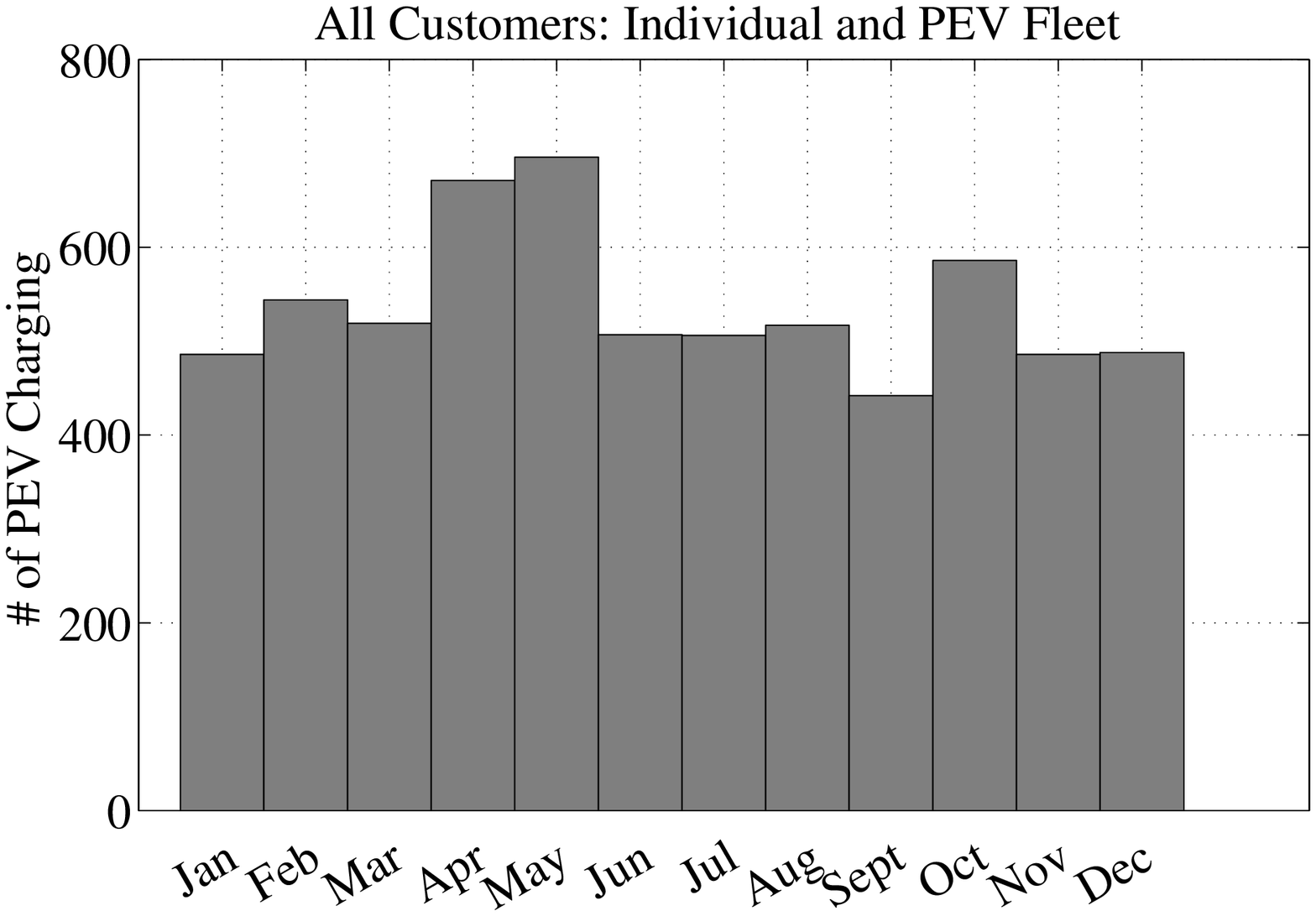}
         \caption{Monthly aggregated demand.}\label{monthBoth}
        \end{subfigure}
        \caption{PEV Monthly demand}\label{monthly}
        \vspace{-5 pt}
\end{figure*}
\begin{figure*}[t]
        \centering
                \begin{subfigure}[b]{0.35\textwidth}
                \centering
                   \includegraphics[width=\columnwidth]{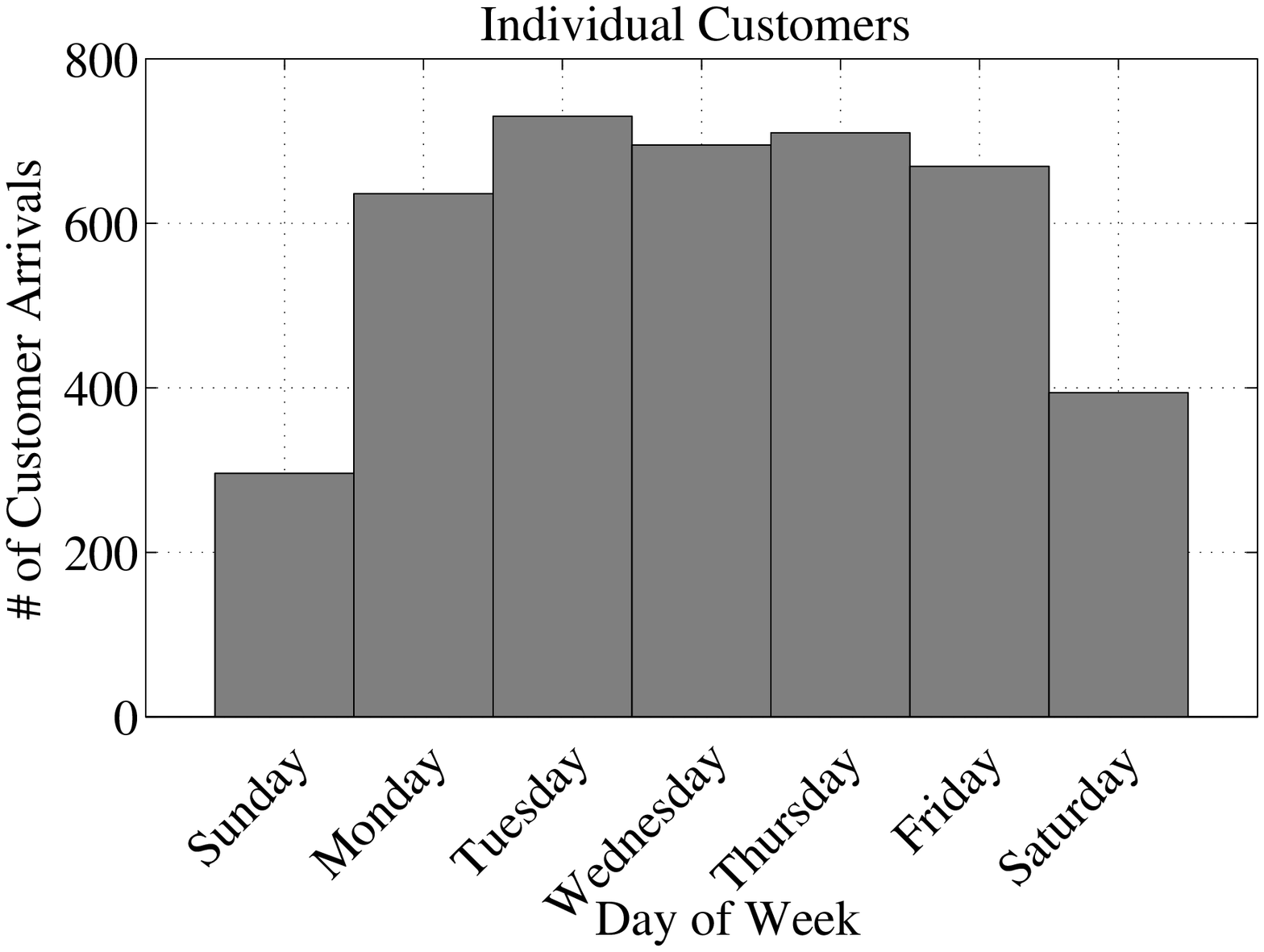}
 \caption{Daily individual customer demand.}\label{dayIndiv}

                      \end{subfigure}
\hspace{-1.57\baselineskip}
        \begin{subfigure}[b]{0.35\textwidth}
                \centering
\includegraphics[width=\columnwidth]{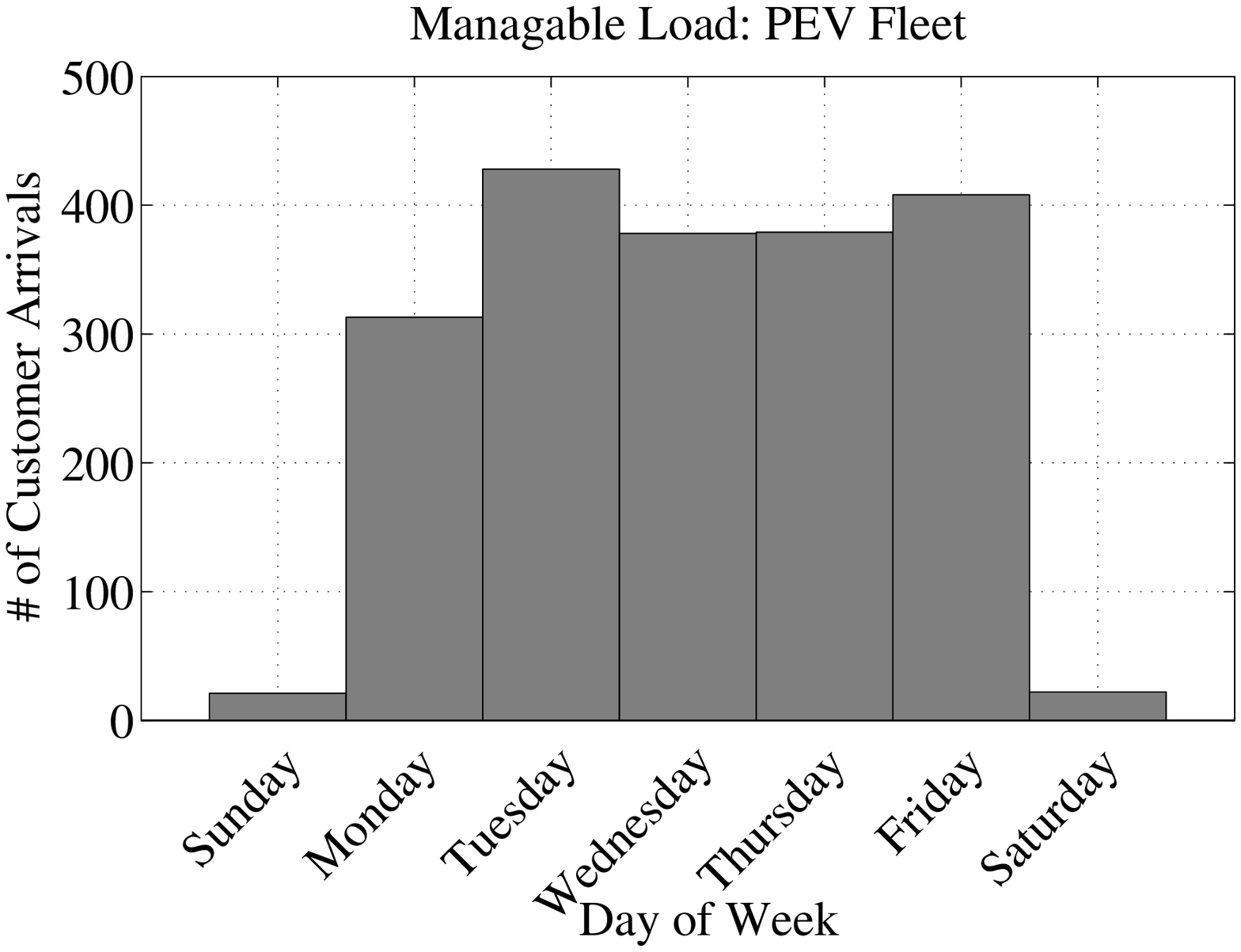}
 \caption{Daily company PEV fleet demand.}\label{dayFleet}

        \end{subfigure}%
        \hspace{-1.45\baselineskip}
        \begin{subfigure}[b]{0.35 \textwidth}
                \centering
                \includegraphics[width=\columnwidth]{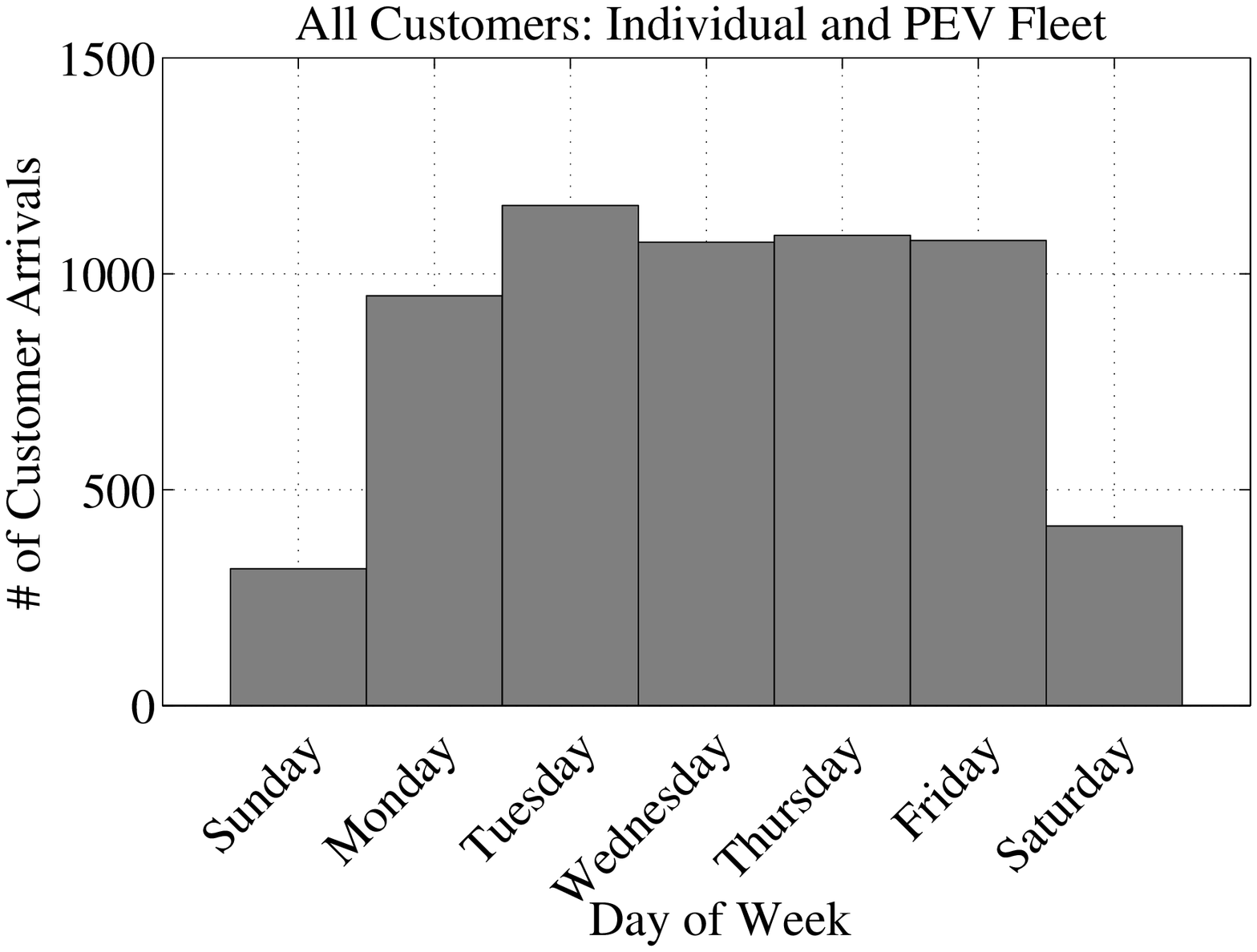}
         \caption{Daily aggregated demand.}\label{dayBoth}
        \end{subfigure}
        \caption{PEV daily demand}\label{daily}
        \vspace{-5 pt}
\end{figure*}
\section{Analysis \& Results}
We start our analysis by differentiating two different PEV types \emph{individuals} and \emph{PEV fleet}. General public such as students, faculty, employees, and visitors constitute the first group, whereas the second group is composed of company vehicles (shuttles, postal cars etc.). The rationale behind this categorization is that since PEV fleets usually have certain tasks, the demand for such vehicles can assumed to be manageable \cite{jsac}. This can be done by optimally adjusting the charging current and service duration, with respect to varying network conditions. Thus, in most of our results are categorized as \emph{individuals}, \emph{fleets}, and \emph{aggregated results}. 

Next, we present the distribution of energy demand of single charging session for each customer case. The results depicted in Fig. \ref{energyDemand} shows that more than half of the customer demand is within $0-8$kWh and $90.1$\% of the total service request is between $0-16$kWh. Considering the fact that mainstream PEV packs range between $16-85$kWh, we conclude that the vast majority of the customers use parking lots to extend their driving ranges. Furthermore, almost $70$\% of the total energy is transferred to individual PEVs. The findings reveal that the system operator can manage the remaining $30$\% of the demand. 

Obviously, these results may not be meaningful without knowing the temporal customer behavior. Thus, in Fig. \ref{hourly} we evaluate the $24$ hours customer behavior. The results in Fig. \ref{hourIndiv} shows that the individual demand peak occurs between $10$am to $12$pm, whereas the demand peak for company fleets is around $1$pm. It is also noteworthy that the peak PEV fleets demand may coincide with the peak demand of other electrical loads. If this trend continues, fleet operators may need to shift their usage to night time. From station design standpoint, this information can be used to decide on how many physical chargers to deploy in order to provide good Level of QoS, as the number of vehicles increases the charging stations may act as \emph{waiting systems}, and station operators could guarantee to finish service within a target deadline.
Moreover, charging facilities are expected to employ solar panels and energy storage units to aid the power grid operations. As given in Section \ref{renewSection}, solar generation depends on daily and seasonal variations. It can be seen that there is a correlation between PEV demand and solar generation (depicted in Fig.\ref{UCSDSolar}). This means that solar energy has a great potential to provide electricity during such hours. Another important parameter is the seasonal customer demand; hence, we show the PEV demand for each month of the year in Fig.\ref{monthly}.
\begin{table}[t]
\caption{Detailed station information.}
\small
\begin{center}
\begin{tabular}{ccccc}
\toprule
 & \# of  & \# of Service  & Dates & Type \\
 & Chargers & Completion & &\\
    \midrule
    St. $1$ & $3$     & $1093$  & $12/19/12$-$1/8/15$ & Public \\
    St. $2$ & $1$     & $196$   & $5/18/13$-$1/7/15$ & Public \\
    St. $3$ & $2$     & $1182$  & $9/14/12$-$1/8/15$ & Public \\
    St. $4$ & $1$     & $838$   & $12/2/11$-$1/8/15$ & Fleet \\
    St. $5$ & $1$     & $147$   & $7/4/13$-$1/8/15$ & Public \\
    St. $6$ & $1$     & $121$   & $5/21/13$-$1/6/15$ & Public \\
    St. $7$ & $2$     & $1294$  & $1/3/13$-$1/8/15$ & Public \\
    St. $8$ & $2$     & $838$   & $12/2/11$-$1/8/15$ & Fleet \\
    St. $9$ & $4$     & $1111$  & $11/15/11$-$1/8/15$ & Public \\
\bottomrule
\end{tabular}
\end{center}
\label{stationInfo1}
\vspace{-20pt}
\end{table}%

Next, we evaluate the daily customer demand. The results depicted in Fig.\ref{daily} shows that the service requests during the weekdays constitute almost a flat profile. From the capacity planning standpoint, this is a desired behavior because less uncertainty in customer demand will translate to higher resource utilization, i.e., storage and solar usage. It is also noteworthy that the customer demand during weekends is considerable less than the weekdays. Station operator can store and sell the energy back to the power grid. Also, some portion of the PEV fleet load on Mondays could be served on Sundays to allow more individual customers.

\begin{figure}[t]
                \centering
                \includegraphics[width=0.75\columnwidth]{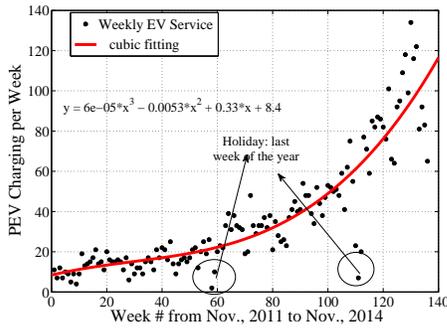}
                 \vspace{-15 pt}

         \caption{Weekly charging data. The number of service request increases as the rise in PEV sale.}\label{weekyData}
        \vspace{-15 pt}
\end{figure}

Our final evaluation is on customer demand over time. We compute the number of PEV services from Nov. $2011$ to Nov. $2013$ and the results depicted in Fig. \ref{weekyData} shows that as the PEV adoption increases the need for charging service surges accordingly. In fact, $4$ out of $9$ charging locations are opened after year $2013$. This data is extremely important to forecast the customer demand and determine the required number of charging nodes. 

\section{Conclusions \& Future Work}
In this paper, we have analyzed and presented the data obtained from a Level $2$ charging network located in major North American University campus. As the future sustainable charging facilities are expected to employ renewable generation and storage units, along with smart energy management systems, we have provided detailed guidelines for charging station design. As a future work, we are aiming to investigate the potential usage of solar powered stations using the obtained data.
\ifCLASSOPTIONcaptionsoff
  \newpage
\fi
  \vspace{-15 pt}
\bibliographystyle{IEEEtran}
\bibliography{references}



\end{document}